\documentclass{article}
\usepackage{amsfonts}
\usepackage{amssymb}
\usepackage{latexsym}
\def\geqs{\geqslant}
\def\leqs{\leqslant}
\def\be{\begin{equation}}
\def\ee{\end{equation}}
\def\bea{\begin{eqnarray}}
\def\eea{\end{eqnarray}}
\def\Qed{\rule{1ex}{1ex}}
\def\ra{\rightarrow}
\def\0{^{\phantom0}}
\def\image{\mathop{\rm im}\nolimits}
\newcommand{\C}{\mathbf{C}}
\newcommand{\R}{\mathbf{R}}
\newcommand{\Z}{\mathbf{Z}}
\newcommand{\PP}{\mathbf{P}}
\def\chim{\chi\0_{H_m}}
\def\hatchim{{\widehat\chi}\0_{H_m}}
\def\muhat{{\widehat\mu}}
\begin{document}
\centerline{\large On finite sequences satisfying linear recursions}
\vspace*{2ex}
\centerline{Noam D.\ Elkies\footnote{
  Supported in part by the Packard Foundation.
  e-mail address: \textsf{elkies@math.harvard.edu}
  }}
\centerline{May, 2001}
  \centerline{\small corrected April, 2002 per referee's report}

\vspace*{5ex}

\begin{quote}
{\bf Abstract.}  For any field $k$\/ and any integers $m,n$
with $0 \leqs 2m \leqs n+1$, let $W_n$ be the \hbox{$k$\/-vector} space
of sequences $(x_0,\ldots,x_n)$, and let $H_m \subseteq W_n$ be the
subset of sequences satisfying a \hbox{degree-$m$} linear recursion,
i.e.\ for which there exist $a_0,\ldots,a_m\in k$\/ such that
subset of sequences satisfying a \hbox{degree-$m$} linear recursion
--- that is, for which there exist $a_0,\ldots,a_m\in k$,
  not all zero, such that
$$
\sum_{i=0}^m a_i x_{i+j} = 0
$$
holds for each $j=0,1,\ldots,n-m$.  Equivalently, $H_m$ is the set
of $(x_0,\ldots,x_n)$ such that the $(m+1) \times (n-m+1)$ matrix
with $(i,j)$ entry $x_{i+j}$ ($0\leqs i\leqs m$, $0\leqs j \leqs n-m$)
has rank at most~$m$.  We use elementary linear and polynomial algebra
to study these sets $H_m$.  In particular, when $k$\/ is a finite field
of~$q$ elements, we write the characteristic function of~$H_m$ as a
linear combination of characteristic functions of linear subspaces
of dimensions~$m$ \hbox{and $m+1$} in~$W_n$.  We deduce a formula
for the discrete Fourier transform of this characteristic function,
and obtain some consequences.  For instance, if the $2m+1$ entries
of a square Hankel matrix of order $m+1$ are chosen independently from
a fixed but not necessarily uniform distribution~$\mu$ on~$k$,
then as $m\ra\infty$ the matrix is singular with probability
approaching $1/q$ provided $\|\muhat\|\0_1 < q^{1/2}$.
This bound $q^{1/2}$ is best possible if $q$ is a square.
\end{quote}

{\large\bf Introduction}

Fix a field $k$.  For any integers $m,n$ with $0 \leqs 2m \leqs n+1$,
let $W_n$ be the \hbox{$k$\/-vector} space of sequences
$(x_0,\ldots,x_n)$, and let $H_m \subseteq W_n$ be the subset
of sequences satisfying a \hbox{degree-$m$} linear recursion,
that is, for which there exist $a_0,\ldots,a_m\in k$, not all zero,
such that
\be
\sum_{i=0}^m a_i x_{i+j} = 0
\label{recursion}
\ee
holds for each $j=0,1,\ldots,n-m$.  Equivalently, $H_m$ is the set
of $(x_0,\ldots,x_n)$ such that the $(m+1) \times (n-m+1)$
Hankel matrix\footnote{
  For more background on Hankel matrices
  (matrices with entries constant on NE-SW diagonals),
  and the closely related Toeplitz or ``persymmetric'' matrices
  (with entries constant on NW-SE diagonals), see for instance
  \cite{Iohvidov}.  These matrices arise in diverse mathematical
  contexts; see for instance \cite{Cantor,DEL} and the references
  in~\cite{Iohvidov}.  In our setting, Hankel matrices are more natural
  than Toeplitz ones, but our results on rank distribution, culminating
  in Thm.~2, apply equally well to matrices of either Hankel or
  Toeplitz type.
  }
\be
\left(
\begin{array}{cccc}
x_0 & x_1 & \cdots & x_{n-m} \\
x_1 & x_2 & \cdots & x_{n-m+1} \\
\vdots & \vdots &  & \vdots \\
x_m & x_{m+1} & \cdots & x_n
\end{array}
\right)
\label{hankmat}
\ee
has rank at most~$m$.\footnote{
  I thank Joe Harris for the geometric observation
  that $H_m$ consists of the lines through the origin coming from
  the points $(x_0:x_1:\cdots:x_n)$ lying on the \hbox{$m$-th}
  secant variety of the rational normal curve
  $(\xi^n : \xi^{n-1} \eta : \cdots : \eta^n)$ in
  \hbox{$n$-dimensional} projective space over~$k$.
  We shall not need this formulation here, but it arises
  naturally in an arithmetic application of~$H_m$~\cite{DEL}.
  }
Now linear recursions on \textit{infinite} sequences
$\{x_i\}_{i\in\Z}$ are known to correspond to polynomials
in the shift operators $T^{\pm 1}: \{x_i\} \mapsto \{x_{i\pm1}\}$,
modulo multiplication by powers of~$T$.  This approach does not work
so nicely for finite sequences, because $T$\/ and $T^{-1}$ push
$x_0$ and $x_n$ off the edge.  We propose to remedy this problem
at $T=0,\infty$ by homogenizing: instead of polynomials in $T^{\pm 1}$,
use homogeneous polynomials in two variables $Y$\/ and~$Z$\/
that act on $W_n$ as the right and left truncation maps to~$W_{n-1}$.
We shall see that this approach yields a clean account of
linear recursions and the subsets $H_m$ in the space $W_n$,
which itself will be identified with the dual of the space $V_n$
of homogeneous polynomials of degree~$n$ in~$Y$\/ and~$Z$.\footnote{
  As an added benefit, the whole structure inherits
  a GL$_2(k)$ structure from the action of GL$_2(k)$
  by linear substitutions on $Y,Z$.
  But this, too, is not needed for the present paper.
  }

In the present paper we develop this account using
elementary linear and polynomial algebra.
When $k$\/ is a finite field of~$q$ elements,
we also write the characteristic function of~$H_m$
as a linear combination of characteristic functions of linear subspaces
of dimensions~$m$ and~$m+1$ in~$W_n$.  We deduce a formula
for the discrete Fourier transform of this characteristic function,
and obtain some consequences.  For instance we obtain a new proof
that $\#H_m = q^{2m}$.  We further show that if the $2m+1$ entries
$x_0,\ldots,x_{2m}$ of a square Hankel matrix
\be
\left(
\begin{array}{cccc}
x_0 & x_1 & \cdots & x_{m} \\
x_1 & x_2 & \cdots & x_{m+1} \\
\vdots & \vdots &  & \vdots \\
x_m & x_{m+1} & \cdots & x_{2m}
\end{array}
\right)
\label{sqhankmat}
\ee
of order $m+1$ are chosen independently from
a fixed but not necessarily uniform distribution~$\mu$ on~$k$,
then as $m\ra\infty$ the the matrix is singular with probability
approaching $1/q$ provided the Fourier transform of~$\mu$ has
$l_1$~norm less than $q^{1/2}$.
This bound is best possible if $q$ is a square:
if $\mu$ is the uniform distribution on $ck_0$,
where $k_0$ is a quadratic subfield of~$k$ and $c\in k^*$ is arbitrary,
then $\|\muhat\|\0_1 = q^{1/2}$ but the probability is $q^{-1/2}$.
It seems reasonable to conjecture that for any $\mu$ the matrix
(\ref{sqhankmat}) is singular with probability $\ra 1/q$ as long as
$\mu$ is supported on a set of at least two elements not contained
in $ck_0$ for any proper subfield $k_0$ of~$k$.

\vspace*{2ex}

{\large\bf The spaces $V_n,W_n$ and some linear algebra}

\textbf{Basic notions and lemmas.}
Fix a field~$k$.  For each integer $n \geq -1$, let $V_n$ be the
vector space of dimension $n+1$ over~$k$ consisting of bivariate
homogeneous polynomials
\be
P(Y,Z) = \sum_{i=0}^n a_i Y^i Z^{n-i}
\label{P_n}
\ee
of degree~$n$.  Let $W_n$ as the dual of~$V_n$.
We identify $W_n$ with the space of sequences $(x_0,\ldots,x_n)$ 
by regarding such a sequence as the linear functional
\be
\sum_{i=0}^n a_i Y^i Z^{n-i} \mapsto \sum_{i=0}^n a_i x_i
\label{x}
\ee
on~$V_n$.  Note that we allow $V_{-1}$ and $W_{-1}$,
each of which is the zero space, but not $V_n,W_n$ for $n < -1$.

If $m \geqs 0$ and $n+1 \geqs m$, polynomial multiplication
$V_m \times V_{n-m} \ra V_n$
gives for each $Q\in V_m$ a linear map $M_n(Q): V_{n-m} \ra V_n$
defined by
\be
M_n(Q): P \mapsto PQ \quad (P \in V_{n-m}).
\label{M_n}
\ee
Our reason for identifying the space $W_n$ of sequences
$(x_0,\ldots,x_n)$ with the dual of $V_n$ is the following observation:

\textbf{Lemma 1.}
\textsl{
  Suppose $Q\in V_m$ is the polynomial $\sum_{i=0}^m a_i Y^i Z^{m-i}$.
  Then the adjoint of~$M_n(Q)$ is the linear map
  $M^*_n(Q): W_n \ra W_{n-m}$ taking any $(x_0,\ldots,x_n)$
  to the sequence of length $n-m+1$ whose \hbox{$j$-th} term is
\be
\sum_{i=0}^m a_i x_{i+j}
\label{M*_n}
\ee
for each $j$ with $0 \leqs j \leqs n-m$.
}

\textit{Proof}\/: We show the equivalent dual statement:
the linear map $M_n(Q)$ takes any polynomial
$$
P(Y,Z) \sum_{j=0}^{n-m} b_j Y^j Z^{n-m-j}
$$
in $V_{n-m}$ to the polynomial $PQ \in V_n$
whose $Y^r Z^{n-r}$ coefficient is $\sum_{i+j=r} a_i b_j$
for each $r$ with $0 \leqs r \leqs n$.
But this is immediate from the expansion of~$PQ$.~~\Qed

Thus $H_m$ is the union of $\ker M^*_n(Q)$ over all nonzero $Q\in V_m$.

Of course that union is not disjoint, but as long as $2m \leqs n+1$
we shall describe the intersection of $\ker M^*_n(Q_1)$ and
$\ker M^*_n(Q_2)$ for any $Q_1,Q_2$ of degree at most~$m$,
see Lemma~4 below.  We first establish some further basic properties:

\textbf{Lemma 2.}
\textsl{
  i) For any $Q,Q'\in V_m$ and $n$ such that $n+1 \geqs m \geqs 0$
we have
\be
M^*_n(Q+Q') = M^*_n(Q) + M^*_n(Q').
\label{Q+Q'}
\ee
 ii) For any $Q_1\in V_{m_1}$, $Q_2\in V_{m_2}$, and $n$ such that
  $m_1,m_2>0$ and $n+1 \leqs m_1+m_2$, we have
\be
M^*_n(Q_1 Q_2)
= M^*_{n-m_2}(Q_1) \circ M^*_n(Q_2)
= M^*_{n-m_1}(Q_2) \circ M^*_n(Q_1).
\label{PQ=QP}
\ee
iii) For any nonzero $Q\in V_m$ and any $n\geq m-1$, the map
$M^*_n(Q)$ is surjective and its kernel has dimension~$m$.
}

\textit{Proof}\/: (i) This is the dual of the identity
$M_n(Q+Q')=M_n(Q)+M_n(Q')$,
which is just the distributive law $P(Q+Q')=PQ+PQ'$
for multiplication of homogeneous polynomials.
(Alternatively, apply Lemma~1.)

(ii) Likewise this is the dual of the fact that multiplying
a polynomial of degree $n-m_1-m_2$ by $Q_1 Q_2$ is the same as
multiplying it first by $Q_2$ and then by $Q_1$ or vice versa.

(iii) Since $k[Y,Z]$ has no zero divisors, $M_n(Q)$ is injective;
thus $M^*_n(Q)$ is surjective, and its kernel has dimension
$\dim W_n - \dim W_{n-m} = m$.~~\Qed

\vspace*{1ex}

\textbf{The ideal $I_x$.}
For $x=(x_0,\ldots,x_n)\in W_n$, define $I_x \subseteq k[Y,Z]$
as follows: any $Q \in k[Y,Z]$ is uniquely $\sum_{m=0}^M Q_m$
with each $Q_m\in V_m$; the subset $I_x$ consists of those
$\sum_{m=0}^M Q_m$ for which $(M^*_n(Q_m))(x)=0$ for each $m \leqs n+1$.

\textbf{Lemma 3.}
\textsl{
  $I_x$ is a homogeneous ideal in $k[Y,Z]$ for all $x\in W_n$.
}

\textit{Proof}\/: By definition $\sum_{m=0}^M Q_m \in I_x$
if and only if each $Q_m \in I_x$.  So it is enough to check
that $I_x \cap V_m$ is closed under addition for each~$m$,
and that $PQ\in I_x$ if $Q\in I_x \cap V_m$ and $P\in V_{m'}$
for some $m,m'\geq 0$.  Each of these is vacuously true if
$m>n+1$ or $m+m'>n+1$ respectively, and follows from 
part~(i) or~(ii) of Lemma~2 otherwise.~~\Qed

The main result of this section is the following partial description
of~$I_x$, stating in effect that it is approximated by a principal ideal
as well as dimension considerations allow:

\textbf{Proposition 1.}
\textsl{
  Suppose for some $x\in W_n$ that $I_x$ contains a nonzero polynomial
of degree at most $(n+1)/2$.
Let $m_0$ be the smallest degree of such a polynomial.  Then
$I_x \cap V_{m_0}$ is {$1$-dimensional}, say
\be
I_x \cap V_{m_0} = k Q_0
\label{Q_0}
\ee
for some nonzero $Q_0 \in V_{m_0}$.  For each $m \leqs n+1-m_0$,
\be
I_x \cap V_m = (M_n(Q_0)) \, (V_{m-m_0}).
\label{IxnVm}
\ee
}%
\textbf{Remark}\/: In particular, it follows that
$I_x \cap V_m$ has dimension $m-m_0+1$ for $m_0 \leqs m \leqs n+1-m_0$.
This cannot hold once $m > n+1-m_0$, except in the trivial case $x=0$,
when $m_0=0$ and $I_x$ is all of $k[Y,Z]$.  Indeed suppose that
$m_0>0$ and $m > n+1-m_0$.  If $m>n+1$ then $I_x \cap V_m = V_m$
has dimension $m+1 > m-m_0+1$.  If $m \leqs n+1$ then $I_x \cap V_m$
is the kernel of the linear map
\be
V_m \ra W_{n-m}, \quad
Q \mapsto (M_n^*(Q))(x);
\label{V->W}
\ee
thus
\be
\dim(I_x \cap V_m) \geqs \dim V_m - \dim W_{n-m} = 2m-n,
\label{2m-n}
\ee
which again exceeds $m-m_0+1$ since $m>n+1-m_0$.  This is what we mean
when we state that $I_x$ approximates the principal ideal $(Q_0)$
as well as dimension considerations allow.

To prove Prop.~1 we must first make good on our promise
to describe intersections of the spaces $\ker M^*_n(Q)$.
We do this in the next lemma, whose statement uses
the greatest common divisor~$Q$
of two homogeneous polynomials $Q_1,Q_2$.
This is defined only up to multiplication by $k^*$,
but such scaling does not affect the space $\ker M^*_n(Q)$,
so the choice of g.c.d.\ will not affect the result.

\textbf{Lemma 4.}
\textsl{
  Let $Q_1,Q_2$ be nonzero polynomials in $V_{m_1},V_{m_2}$
respectively, with greatest common divisor~$Q$.  Then, for each
$n \geqs \max(m_1,m_2) - 1$,
\be
\ker M^*_n(Q_1) \cap \ker M^*_n(Q_2) \supseteq \ker M^*_n(Q),
\label{lemma4}
\ee
with equality if and only if
\be
n+1 \geqs m_1 + m_2 - \deg(Q).
\label{lemma4eq}
\ee
}%
\textit{Proof}\/: If $x \in \ker M^*_n(Q)$ then $x$ is in the kernel
of both $M^*_n(Q_1)$ and $M^*_n(Q_2)$, because each of these linear maps
factors through $M^*_n(Q)$ by part~(ii) of Lemma~2.
Thus $x$ is in the intersection of the two kernels,
whence (\ref{lemma4}) follows.
It remains to establish the condition of equality.

Let $m=\deg Q$, and $m'=m_1+m_2-m$.  By Lemma 2(iii),
the codimensions in~$W_n$ of $\ker M^*_n(Q_1)$ and $\ker M^*_n(Q_2)$
are $n+1-m_1$ and $n+1-m_2$ respectively.  Thus their intersection
has codimension at most
\be
(n+1-m_1) + (n+1-m_2) = (n+1-m) + (n+1-m').
\label{codim_bound}
\ee
Hence if $m'>n+1$ then this codimension is strictly less than
the codimension of $\ker M^*_n(Q)$.  Thus the condition $m' \leqs n+1$
is necessary for equality in~(\ref{lemma4}).

We conclude the proof by showing that this condition is also sufficient.
Let $Q'$ be the least common multiple
\be
Q' = Q_1 Q_2 / Q
\label{Q'}
\ee
of $Q_1$ and $Q_2$; this is a homogeneous polynomial of degree~$m'$.
Assuming that $m' \leqs n+1$, we may then consider $M^*_n(Q')$.
We claim that
\be
\ker M^*_n(Q_1) + \ker M^*_n(Q_2) = \ker M^*_n(Q').
\label{lemma4'}
\ee
By duality, this claim is equivalent to
\be
\image(M_n(Q_1)) \cap \image(M_n(Q_2)) = \image(M_n(Q')).
\label{lemma4'dual}
\ee
But this is just the statement that a polynomial in~$V_n$
is divisible by both $Q_1$ and $Q_2$
if and only if it is divisible by $Q'$ ---
which is true because $Q'$ is the least common multiple
of~$Q_1$ and~$Q_2$.  We thus have
\bea
&& \dim (\ker M^*_n(Q_1) \cap \ker M^*_n(Q_2))
\nonumber \\
&\!\!\!\!\!=& \dim (\ker M^*_n(Q_1)) + \dim (\ker M^*_n(Q_2))
 - \dim (\ker M^*_n(Q_1) + \ker M^*_n(Q_2))
\nonumber \\
&\!\!\!\!\!=& \dim (\ker M^*_n(Q_1)) + \dim (\ker M^*_n(Q_2))
 - \dim (\ker M^*_n(Q')).
\label{dosido}
\eea
By Lemma 2(iii) again, this dimension equals
\be
m_1 + m_2 - m' = m = \dim(\ker M^*_n(Q)).
\label{aha4}
\ee
Since we already know that $\ker M^*_n(Q_1) \cap \ker M^*_n(Q_2)$
contains $\ker M^*_n(Q)$, we conclude that these two spaces
are equal.~~\Qed

\textbf{Corollary.}
\textsl{
  Suppose $x\in W_n$.  If $I_x$ contains homogeneous polynomials
$Q_1,Q_2$ whose least common multiple has degree at most $n+1$,
then $I_x$ contains $\gcd(Q_1,Q_2)$.  In particular, this conclusion
holds if $\deg Q_1 + \deg Q_2 \leqs n+1$.
}

\textit{Proof}\/: Under our hypotheses, $x$ is contained in
both $\ker M^*_n(Q_1)$ and $\ker M^*_n(Q_2)$,
and the equality condition of Lemma~4 is satisfied.
Therefore
\be
x \in \ker M^*_n(Q_1) \cap \ker M^*_n(Q_2) = \ker M^*_n(\gcd(Q_1,Q_2)),
\label{gcd}
\ee
which is to say that $I_x$ contains $\gcd(Q_1,Q_2)$ as claimed.~~\Qed

We can now easily prove Prop.~1.  Suppose $Q_1,Q_2$ are nonzero
polynomials in $I_x \cap V_{m_0}$.  By the hypothesis of Prop.~1
we know $2m_0\leqs n+1$.  The Corollary to Lemma~4 thus applies,
and we find that $I_x$ contains $\gcd(Q_1,Q_2)$.  Unless $Q_1,Q_2$
are proportional, $\deg(\gcd(Q_1,Q_2))<m_0$, which is impossible
by the definition of~$m_0$.  Thus $I_x \cap V_{m_0}$ has dimension~$1$
as claimed.  By the same Corollary, if $m\leqs n+1-m_0$ and
$Q\in I_x \cap V_{m_0} - \{0\}$ then $I_x \ni \gcd(Q_0,Q)$.
Since again $\gcd(Q_0,Q)$ must have degree at least~$m_0$,
we conclude that $Q$\/ is a multiple of~$Q_0$.  Since $I_x$
is an ideal (Lemma~3), we already know that $I_x$ contains
all multiples of~$Q_0$; thus $I_x \cap V_{m_0}$ consists of
all \hbox{degree-$m$} multiples of~$Q_0$, and we are done.~~\Qed\,\Qed

It is thus natural to call $Q_0$ the \textsl{minimal linear recursion}
satisfied by~$x$.  (Again $Q_0$ is defined only up to multiplication
by~$k^*$.) From Prop.~1 we deduce the following description
of the degree~$m_0$ of this minimal recursion:

\textbf{Corollary.}
\textsl{
  If $x\in H_m$ for some $m \leqs (n+1)/2$ then the degree of
the minimal linear recursion satisfied by~$x$ equals the rank
of the Hankel matrix~(\ref{hankmat}) associated to~$x$.
}

\textit{Proof}\/: Let $m_0$ be this minimal degree.
The rank of (\ref{hankmat}) is $m+1-d$, where
$d$\/ is the dimension of the kernel of the action of this matrix
on row vectors of length~$m+1$.  But Lemma~1 identifies this kernel
with the space $I_x \cap V_m$ of \hbox{degree-$m$} recursions
satisfied by~$x$.  Since $m_0 \leqs m \leqs (n+1)/2$,
we may apply Prop.~1 to find that $d = m - m_0 + 1$.
Thus $m_0$ is the rank of the Hankel matrix, as claimed.~~$\Qed$

\vspace*{2ex}

{\large\bf The characteristic function of $H_m$}

\textbf{Decomposition into signed linear subspaces.}
We assume henceforth that $k$ is a finite field of $q$ elements.
For integers $m,n$ satisfying our customary condition $2m \leqs n+1$,
let $\PP_m$ be the set of all subspaces of~$W_n$ of the form
$\ker M^*_n(Q)$ for some nonzero $Q \in V_m$.
(By Lemma~4 and part~(iii) of Lemma~2,
 $\ker M^*_n(Q_1) = \ker M^*_n(Q_2)$ if and only if
$Q_1,Q_2$ are proportional; thus $\PP_m$ consists of
\be
\frac{\#(V_m - \{0\})}{\#(k^*)} = \frac{q^{m+1}-1}{q-1}
\label{|PP_m|}
\ee
subspaces.  We note for later use that this formula remains valid
if we allow $m=-1$, when $\PP_m$ is empty.)
Recall that we defined $H_m$ as the set of $x\in W_n$
satisfying a recursion of degree~$m$,
and noted that $H_m$ is thus the union of all the subspaces in~$\PP_m$.
We further noted that this union is not disjoint,
and thus that $\chim$, the characteristic function of~$H_m$,
is not simply the sum of the characteristic functions
of the subspaces in $\PP_m$.  However, by Lemma~4,
the intersection of any two subspaces in $\PP_m$ is again
the kernel of $M^*_n(Q)$ for some nonzero homogeneous~$Q$\/
of degree $\leqs m$, and more generally if $m_1,m_2\leqs m$ then 
the intersection of any subspace in $\PP_{m_1}$
with any subspace in $\PP_{m_2}$ is itself in $\PP_{m'}$
for some $m' \leqs m$.  Thus we can use inclusion-exclusion identities
to write $\chim$ as a linear combination of the characteristic functions
of subspaces in $\PP_{m'}$ for $m'\leqs m$.
Fortunately the resulting formula is quite simple:

\textbf{Proposition 2.}
\textsl{
  The characteristic function of~$H_m$ equals
\be
\sum_{K \in \PP_m} \chi_K\0
\; - \; q \!\!\sum_{K \in \PP_{m-1}} \chi_K\0,
\label{prop2}
\ee
in which $\chi_K\0$ is the characteristic function of the set~$K$,
and the second sum is interpreted as zero when $m=0$.
}

\textit{Proof}\/: Clearly (\ref{prop2}) is an integer-valued function
on~$W_n$ supported on~$H_m$.  Thus we need only show that its value
at~$x$ equals~$1$ for all $x\in H_m$.  But this value is
\bea
&& \frac{\#(I_x \cap V_m) - 1}{q-1}
\, - \, q \, \frac{\#(I_x \cap V_{m-1}) - 1}{q-1}
\nonumber \\
&=& 1 + \frac{\#(I_x \cap V_m) - q\, \#(I_x \cap V_{m-1})} {q-1} \ .
\label{f(x)}
\eea
Let $m_0$ be the degree of the minimal linear recursion satisfied
by~$x$.  By Prop.~1, $I_x \cap V_m$ and $I_x \cap V_{m-1}$ are
vector spaces of dimensions $m-m_0+1$ and $m-m_0$ respectively over~$k$.
(Note that this remains true if $m_0=m$, when $I_x \cap V_{m-1}$
is the zero space.)  Thus $\#(I_x \cap V_m) = q\; \#(I_x \cap V_{m-1})$,
and (\ref{f(x)}) simplifies to~$1$ as claimed.~~\Qed\,\Qed

We easily deduce the formula~\cite[Thm.~1]{Daykin}
for the size of~$H_m$:

\textbf{Corollary.}
\textsl{
  For all nonnegative $m \leqs (n+1)/2$ we have
\be
\#(H_m) = q^{2m}.
\label{|H_m|}
\ee
}%
\textit{Proof}\/:
The size of $H_m$ is the sum of~$\chim(x)$ over $x\in W_n$.
By (\ref{prop2}), this sum is
\be
\sum_{K \in \PP_m} \!\#(K)
\; - \; q \!\!\sum_{K \in \PP_{m-1}} \!\#(K).
\label{prop2cor}
\ee
But by Lemma 2(iii), each $K\in\PP_m$ has size $q^m$,
and each $K\in\PP_{m-1}$ has size $q^{m-1}$.
Using (\ref{|PP_m|}) --- and this is where we use the validity
of (\ref{|PP_m|}) \hbox{also for $m=-1$} ---
we thus simplify (\ref{prop2cor}) to
\be
\frac{q^{m+1}-1}{q-1} \, q^m \; - \; q \, \frac{q^m-1}{q-1} \, q^{m-1}
= q^{2m},
\label{prop2pf}
\ee
as claimed.~~\Qed

In particular, if $n=2m-1$ then $\#H_m = \#W_n$, whence $H_m = W_n$ ---
which is clear because in this case the Hankel matrix (\ref{hankmat})
has only $m$ rows, so must have rank at most~$m$.
(This is essentially the special case $n=2m-1$ of the dimension count
we used earlier to deduce (\ref{2m-n}); in this case we find that
$I_x \cap V_m$ has rank at least $2m-n=1$,
so must contain a nonzero vector.)
Starting from this, one may establish without too much difficulty
a bijection from $W_{2m-1}$ to the subset $H_m$ of~$W_n$
for any $n \geqs 2m-1$, even without our $k[Y,Z]$ framework.
(This is in effect how (\ref{|H_m|}) is proved in~\cite{Daykin}.)
But our approach also yields a formula for the Fourier transform
$\hatchim(P)$ for all $P\in V_n$, whereas
(\ref{|H_m|}) only gives $\hatchim(0)$.  We turn to $\hatchim$ next.

\vspace*{1ex}

\textbf{Discrete Fourier transform.}
To define the Fourier transform on~$W_n$, we first define it on~$k$.
Fix a nontrivial character~$\psi_0$ of~$k$, that is,
a nontrivial homomorphism from the additive group of~$k$\/
to the unit circle in~$\C$.  [If $k = \Z/p\Z$ for some prime~$p$,
we may take $\psi_0(x) = \exp (2\pi i x / p)$; in general $k$\/
contains $\Z/p\Z$ where $p$ is the characteristic of~$k$,
and we may take $\psi_0(x) = \exp (2\pi i t(x) / p)$ where
$t : k \ra \Z/p\Z$ is any nontrivial homomorphism of additive groups.
One common choice for~$t$\/ is the trace from~$k$\/ to~$\Z/p\Z$.
At any rate none of our results will depend on the choice of~$\psi_0$.]
For any function $f: k \ra \C$, we define
the (discrete) \textit{Fourier transform} $\widehat f$\/ of~$f$\/
to be the following function from~$k$\/ to~$\C$\/:
\be
\widehat f(a) := \sum_{x \in k} f(x) \psi_0(ax).
\label{f2fhat}
\ee
It is known that $f \mapsto \widehat f$\/ is a linear bijection
on the space~$\C^q$ of complex-valued functions on~$k$, and that
the inverse bijection is given by the
\textit{Fourier inversion formula}:
\be
f(x) = \frac1q \sum_{a \in k} \widehat f(a) \psi_0(-ax).
\label{fhat2f}
\ee
The Fourier transform is defined more generally for finite-dimensional
vector spaces over~$k$.  Let $V,W$\/ be a dual pair of such spaces,
of dimension~$d$.  (We shall use $V=V_n$, $W=W_n$, $d=n+1$.)
To each function $F: W \ra \C$\/ we associate
its discrete Fourier transform
\be
\widehat F(a) := \sum_{x \in W} F(x) \psi_0(\langle a,x \rangle).
\label{F2Fhat}
\ee
Again $F \mapsto \widehat F$\/ is a linear bijection,
and in this context the inversion formula reads
\be
F(x) = \frac1{q^d}
 \sum_{a \in V} \widehat F(a) \psi_0(-\langle a,x \rangle).
\label{Fhat2F}
\ee
To recover $\hatchim$ from Prop.~2, we shall need one more fact
about the discrete Fourier transform:

\textbf{Lemma 5.}
\textsl{
  For any linear subspace $K \subseteq W$,
the Fourier transform of its characteristic function~$\chi\0_K$
is $(\#K) \cdot \chi\0_{K^\perp}$,
where $K^\perp$ is the annihilator of~$K$\/ in~$V$.
}

\textsl{Proof}\/:
By definition, ${\widehat\chi}\0_K(a)$ is the sum over~$K$\/
of the character $x \mapsto \psi_0(\langle a,x \rangle)$;
thus ${\widehat\chi}\0_K(y) = \#K$ or~$0$
according as this character is trivial or nontrivial on~$K$,
that is, according as $a \in K^\perp$ or $a \notin K^\perp$.~~\Qed

We can now give our formula for $\hatchim$.
It will be convenient to introduce the following notation:
for $P\in V_n$ and any integer~$d$, define $\omega_d(P)$
to be $1/(q-1)$ times the number of nonzero $Q\in V_d$
such that $P$\/ is a multiple of~$Q$.  Equivalently,
$\omega_d(P)$ is the number of \hbox{degree-$d$}\/ factors of~$P$\/
up to $k^*$ scaling, and the number of homogeneous principal ideals
in $k[Y,Z]$ that contain~$P$\/ and have a generator of degree~$d$.
For instance, $\omega_0(P)=1$, and for all $d \geq -1$,
\be
\omega_d(0) = \frac{q^{d+1}-1}{q-1}
\; [ = \#(\PP_d) \ {\rm if}\ 2d\leqs n+1].
\label{om(0)}
\ee
Moreover, for nonzero~$P$\/ we have the identity
\be
\omega_d(P) = \omega_{n-d}(P),
\label{d,n-d}
\ee
due to the bijection $Q \leftrightarrow P/Q$\/
between factors of~$P$\/ of degree~$d$\/ and~$n-d$.
(The notation $\omega_d$ is suggested by the omega function
in elementary number theory, which counts the positive divisors
of a given positive integer.)

\textbf{Theorem 1.}
\textsl{
  For every $m \leqs (n+1)/2$ and $P\in V_n$ we have
\be
\hatchim(P) = q^m \left( \omega_m(P) - \omega_{m-1}(P) \right).
\label{hatchim}
\ee
}%
\textsl{Proof}\/:
By Prop.~2 and Lemma~5, this follows from the following observation:
for any homogeneous polynomial~$Q$\/ of degree at most~$n$,
the annihilator in~$V_n$ of $\ker M^*_n(Q)$ is the image of
$M_n(Q)$, which is the space of \hbox{degree-$n$} multiples of~$Q$.
Thus when we use (\ref{prop2}) to expand $\hatchim$ as a linear
combination of characteristic functions of annihilators,
the number of subspaces in $\PP_m$ or $\PP_{m-1}$
that contribute a term to $\hatchim(P)$
is the number of divisors of~$P$\/ of degree~$m$ or~$m-1$
up to $k^*$ scaling.  Each of these terms is $q^m$ or
$-q\cdot q^{m-1} = -q^m$ respectively,
whence the formula~(\ref{hatchim}).~~\Qed\,\Qed

As promised, Prop.~2 is the special case $P=0$ of this formula
(cf.~(\ref{om(0)})).  Also, if $n=2m-1$, the identity (\ref{d,n-d})
yields $\hatchim(P)=0$ for all $P\neq0$, consistent with $H_m=W_n$
in that case.

\vspace*{1ex}

\textbf{Hankel matrices with independently biased entries.}
The formula~(\ref{|H_m|}) can be interpreted thus:
if $x_0,\ldots,x_n$ are chosen independently at random
from the uniform distribution on~$k$, then the resulting vector
$(x_0,\ldots,x_n)$ is in~$H_m$ with probability $q^{2m-(n+1)}$.
Using Thm.~1 we can also get at the probability that
$(x_0,\ldots,x_n) \in H_m$ if the $x_i$ are still chosen
independently at random but from distributions $\mu_i$ on~$k$\/
that are not necessarily uniform.

We regard the $\mu_i$ as functions from~$k$ to~$\R$\/
satisfying the conditions:
$\mu_i(x) \geqs 0$ for all $x\in k$, and
\be
[\muhat_i(0) = ]\; \sum_{x\in k} \mu_i(x) = 1.
\label{sum=1}
\ee
Then the probability that $\vec x := (x_0,\ldots,x_n)$ is in~$H_m$ is
\be
\Pi_m(\mu_0,\ldots,\mu_n) =
\sum_{\vec x \in W_n} \chim(\vec x) \prod_{i=0}^n \mu_i(x_i).
\label{muprob}
\ee
By applying Fourier inversion to~$\chim$ we can express this
as a linear combination of the values of~$\hatchim(P)$.
The resulting formula is:

\textbf{Lemma 6.}
\textsl{
  We have
\be
\Pi_m(\mu_0,\ldots,\mu_n) =
q^{-(n+1)} \sum_{P \in V_n} \hatchim(P) \prod_{i=0}^n \muhat_i(-a_i),
\label{muprob6}
\ee
where $a_i$ is the $Y^i Z^{n-i}$ coefficient of~$P$\/ as in~(\ref{P_n}).
}

\textit{Proof}\/: By Fourier inversion~(\ref{Fhat2F}),
\be
\Pi_m(\mu_0,\ldots,\mu_n) =
q^{-(n+1)} \sum_{P \in V_n} \hatchim(P)
\left( \sum_{\vec x \in W_n}
  \psi_0(-\langle P,x \rangle) \prod_{i=0}^n \mu_i(x_i)
\right).
\label{invert}
\ee
Now $\langle P,x \rangle = \sum_{i=0}^n a_i x_i$, so
\be
\psi_0(-\langle P,x \rangle) \prod_{i=0}^n \mu_i(x_i)
= \prod_{i=0}^n \psi_0(-a_i x_i) \mu_i(x_i).
\label{psi=prod}
\ee
Thus the inner sum in~(\ref{invert}) factors into
\be
\prod_{i=0}^n \left( \sum_{x_i\in k} \psi_0(-a_i x_i) \mu_i(x_i) \right)
= \prod_{i=0}^n \muhat_i(-a_i).
\label{prodmu}
\ee
Entering this into~(\ref{invert}) yields the claimed
formula~(\ref{muprob6}).~~\Qed

The term $P=0$ in~(\ref{invert}) contributes
\be
q^{-(n+1)} \hatchim(0) \prod_{i=0}^n \muhat_i(0) = q^{2m-(n+1)},
\label{hatchim(0)}
\ee
because $\hatchim(0)=q^{2m}$ and each $\muhat_i(0)=1$.
The absolute value of the sum of the remaining terms is at most
\bea
&& q^{-n+1} \sup_{P\in V_n - \{0\}} |\hatchim(P)|
\ \cdot \sum_{P\in V_n - \{0\}} \, \prod_{i=0}^n |\muhat_i(-a_i)|
\nonumber \\
&=& q^{-n+1} \sup_{P\in V_n - \{0\}} |\hatchim(P)|
\left[ \left( \prod_{i=0}^n \| \muhat_i \|\0_1 \right) - 1 \right],
\label{prodmuhat}
\eea
where $\| \muhat_i \|\0_1$ is the $l_1$ norm
\be
\| \muhat_i \|\0_1 := \sum_{a\in k} |\muhat_i(a)|.
\label{l1}
\ee
Since $\muhat_i(0)=1$, we have $\| \muhat_i \|\0_1 \geqs 1$,
with equality if and only if $\muhat_i(a)=0$ for all $a\neq 0$.
By Fourier inversion~(\ref{fhat2f}), this condition is equivalent
to $\mu_i(x)=1/q$ for all~$x$.  Hence $\| \muhat_i \|\0_1 = 1$
if and ony if $\mu_i$ is the uniform distribution on~$k$.
We may thus regard $(\prod_{i=0}^n \| \muhat_i \|\0_1) - 1$
as a measure of how far the product distribution $\mu_0 \cdots \mu_n$
departs from uniform distribution on~$W_n$.

What of the other factor $\sup_{P\neq0} |\hatchim(P)|$
in the error estimate (\ref{prodmuhat})?  By Thm.~1, each $\hatchim(P)$
is a multiple of~$q^m$.  Once $n \geqs 2m$, we cannot expect
$\hatchim(P)$ to vanish for all $P\neq 0$, so
$\sup_{P\neq0} |\hatchim(P)|$ must be at least $q^m$.
We next show that it $|\hatchim(P)|$ is never much larger than $q^m$
for $P\neq 0$:

\textbf{Lemma 7.}
\textsl{
  For every $q$ and $\epsilon>0$,
there exists an effective constant~$C$\/ such that
\be
\omega_d(P) < C (1+\epsilon)^n
\label{omega_bd}
\ee
for every nonzero $P\in V_n$ and every integer~$d$.
}

(This is analogous to the standard fact that the number of factors
of an \hbox{$n$-digit} integer is subexponential in~$n$,
and will be proved in the same way.)

\textit{Proof}\/: Define
\be
\omega(P) := \sum_{d=0}^n \omega_d(P),
\label{omega}
\ee
the total number of divisors of~$P$\/ up to~$k^*$ scaling.
Factor $P$\/ into irreducibles over~$k$\/:
\be
P = \prod_{s=1}^r P_s^{e_s},
\label{factor}
\ee
with $P_s$ distinct irreducibles of degree~$f_s$.
Comparing degrees in~(\ref{factor}) we find
\be
n = \sum_{s=1}^r e_s f_s.
\label{sumef}
\ee
Now
\be
\omega(P) = \prod_{s=1}^r (e_s+1),
\label{omegaprod}
\ee
because the general divisor of~$P$\/ is $\prod_{s=1}^r P_s^{e'_s}$
with each $e'_s$ chosen from among the $e_s+1$ possibilities
$0,1,\ldots,e_s$.
Fix $m_0$ large enough that $2^{1/m_0} < 1+\epsilon$,
and factor (\ref{omegaprod}) as
\be
\omega(P) = \prod_{f_s < m_0} (e_s+1) \prod_{f_s \geqs m_0} (e_s+1).
\label{omegaprod2}
\ee
The second product is at most
\be
\prod_{f_s \geqs m_0} 2^{e_s}
= 2^{\sum_{f_s \geqs m_0} e_s}_{\phantom0}
\leqs 2^{n/m_0},
\label{prod.2}
\ee
since $m_0 \sum_{f_s \geqs m_0} e_s \leqs \sum_{s=1}^r e_s f_s = n$
by~(\ref{sumef}).  The first product in (\ref{omegaprod2})
has at most $B$\/ factors, where $B$\/ is the number of
irreducible bivariate homogeneous polynomials of degree $<m_0$
up to $k^*$~scaling.  Each factor is at most $n+1$, so
the product is at most $(n+1)^B$.  Since $\log (n+1)^B = o(n)$
as $n\ra\infty$, and $2^{1/m_0} < 1+\epsilon$,
we conclude that
\be
\omega(P) \leqs 2^{n/m_0} (n+1)^B \ll (1+\epsilon)^n.
\label{lemma7pf}
\ee
Since $\omega_d(P) \leqs \omega(P)$,
we deduce $\omega_d(P) \ll (1+\epsilon)^n$.~~\Qed

Combining this estimate with Thm.~1 and Lemma~6, we obtain:

\textbf{Theorem 2.}
\textsl{
  For every $q$ and $\epsilon>0$,
there exists an effective constant~$C$\/ such that
\be
\left| \Pi_m(\mu_0,\ldots,\mu_n) - q^{2m-(n+1)} \right|
< C (1+\epsilon)^n q^{m-n} \prod_{i=0}^n \| \muhat_i \|\0_1 \, .
\label{thm2}
\ee
for any $n$ and any distributions~$\mu_i$ on~$k$.~~\Qed\,\Qed
}

In particular, suppose that $n=2m+\alpha$ for some fixed nonnegative
integer~$\alpha$, and that all the $\mu_i$ are the same,
so that each $x_i$ is chosen from the same distribution~$\mu$.
Then, as long as $\| \muhat \|\0_1 < q^{1/2}$,
the error term in~(\ref{thm2}) approaches~$0$ as $m\ra\infty$,
and we conclude that if each of $x_0,\ldots,x_{2m+\alpha}$
is chosen independently from the distribution~$\mu$ then
$x \in H_m$ with probability approaching $q^{-(\alpha+1)}$,
same as for the uniform distribution.  As noted in the Introduction,
the bound on $\| \muhat \|\0_1$ is best possible,
at least if $q$ is a square:
in that case $k$ has a quadratic subfield~$k_0$,
and if each $x_i$ is chosen uniformly from $k_0$
(or from $ck_0$ for some $c\in k^*$) then $x\in H_m$
with probability $q^{-(\alpha+1)/2}$, not $q^{-(\alpha+1)}$;
but for this distribution, $\| \muhat \|\0_1 = q^{1/2}$ by Lemma~5.


\pagebreak

{\large\bf Open questions}

\textbf{Better bounds on $\Pi_m(\mu_0,\ldots,\mu_n) - q^{2m-(n+1)}\;$?}
We showed (Thm.~2) that $\Pi_m(\mu_0,\ldots,\mu_n)$ is well approximated
by $q^{2m-(n+1)}$ under certain hypotheses on the $\mu_i$.
Can these hypotheses by weakened by lowering the error bound
in~(\ref{thm2})?  Of course we must exclude some choices of $\mu_i$.
For instance we certainly cannot have every $\mu_i$ supported
on only one point; and we already gave the counterexample
of uniform distribution on a proper subfield of~$k$.
But it seems plausible that, except for such pathological cases,
$(x_0,\ldots,x_n)$ should be about as likely to be in~$H_m$
with $x_i$ chosen from~$\mu_i$ as it is with $\vec x$ chosen
uniformly from~$W_n$ --- whether or not the $\| \mu_i \|\0_1$ are
small enough to deduce $\Pi_m(\mu_0,\ldots,\mu_n) \sim q^{2m-(n+1)}$
from Thm.~2.  For instance we may surmise the following

\textbf{Conjecture.}
\textsl{
  Fix $k$\/ and a closed set~$K$\/ of distributions $\mu : k\ra\R$.
Assume that no $\mu\in K$\/ is supported on a single point, nor 
on $ck_0$ for any $c\in k^*$ and any proper subfield~$k_0$ of~$k$.
Then, for every real $R\geqs 2$, we have
\be
\Pi_m(\mu_0,\ldots,\mu_n) = (1+o(1)) q^{2m-(n+1)}
\label{conj}
\ee
for any sequence of $(n,m,\mu_0,\ldots,\mu_n)$ for which
$m \ra \infty$,  $2m \leqs n \leqs Rm$, and $\mu_i\in K$\/ for each~$i$.
}

In particular, suppose $q=R=2$.  A distribution on~$k$\/ is then
a pair $(\mu(0),\mu(1))$ of nonnegative numbers with $\mu(0)+\mu(1)=1$.
The conjecture then asserts that, for each $p>0$, if each entry~$x_i$
of a square Hankel matrix of order~$m+1$ over $\Z/2\Z$ is chosen
independently at random with probabilities $\mu_i(0),\mu_i(1)$
both $\geqs p$, then the matrix is singular
with probability approaching~$1/2$ as $m\ra\infty$.
Thm.~2 shows this only for $p > 1 - 2^{-1/2} \approx 29.3\%$.

\vspace*{1ex}

\textbf{Higher dimensions.}
What happens to our theory
in the context of arrays of dimension~$2$ or greater,
rather than finite sequences?
One could start the analysis in the same way, using for instance
homogeneous polynomials in three variables to treat triangular arrays,
or bihomogeneous polynomials in two pairs of variables for
rectangular arrays.  The resulting structures
will surely be more complicated in higher dimensions,
but it may still be possible to find tractable descriptions.

\vspace*{1ex}

\textbf{Determinants of nonsingular Hankel matrices.}
In another direction, we return to the case $n=2m$
of square Hankel matrices~(\ref{sqhankmat}) of order~$m+1$,
for which $H_m$ consists in effect of such matrices
whose determinant vanishes.
We then ask: is there a formula analogous to~(\ref{hatchim}),
or even an estimate analogous to Lemma~7,
for the discrete Fourier transform of the set of square Hankel matrices
of order~$m+1$ with determinant~$c$,
for any given {\em nonzero} $c\in k$\/?
This is easy when $q=2$,
in which case that set is just the complement of~$H_m$.
But the problem seems to require new techniques once $q\geqs3$.

\end{document}